\documentclass[10pt]{article}
\usepackage{}
\usepackage{cite}
\usepackage{mathrsfs}
\usepackage{amsfonts}
\usepackage{amsmath}
\usepackage{amsfonts,amssymb,color}
\usepackage{dsfont}
\usepackage{curves}
\usepackage{mathrsfs}
\usepackage{pifont}
\usepackage{amssymb}
\usepackage{graphicx}
\usepackage{float}

\numberwithin{equation}{section}
\newtheorem{theorem}{Theorem}[section]

\newtheorem{definition}{Definition}[section]
\numberwithin{figure}{section}
\newtheorem{lemma}{Lemma}[section]
\newtheorem{corollary}{Corollary}[section]
\newtheorem{remark}{Remark}[section]

\newcommand{\qed}{\hfill\rule{0.5em}{0.809em}}

\def\emptyset{\mbox{{\rm \O}}}

\def\qed{\hfill \rule{4pt}{7pt}}

\begin{document}

\title{Extremal graphs for  vertex-degree-based  invariants \\ with given degree  sequences  \thanks{ The first author is supported by   NSFC of China (No. 11571123), the Training Program for
Outstanding Young Teachers in University of Guangdong Province (No. YQ2015027),   Guangdong Engineering Research Center for Data Science (No. 2017A-KF02) and Guangdong Province Ordinary University Characteristic Innovation Project
(No. 2017KTSCX020).  The second author is supported by NNSF of China (No. 11671202) and Chinese Excellent Overseas Researcher Funding in 2016. The third  author is supported by    the Joint NSFC-ISF Research Program (jointly funded by the National Natural Science Foundation of China and the Israel Science Foundation (No. 11561141001)) and  the National Natural Science Foundation of China (No. 11531001).  \newline  E-mail
addresses: liumuhuo@163.com (M. Liu),   kexxu1221@126.com(K. Xu),  xiaodong@sjtu.edu.cn (X.-D. Zhang, Corresponding author)}}
  \author{\small Muhuo Liu$^{1,2}$,  Kexiang Xu$^{3}$, Xiao-Dong Zhang$^{4}$,
\\\small $^{1}$  Department of Mathematics, College of Mathematics and Informatics, \\ \small South China Agricultural University, Guangzhou, 510642, China \\\small $^{2}$  College of Mathematics and Statistics, Shenzhen University, \\\small Shenzhen, 518060,  China
\\\small $^{3}$ College of Science, Nanjing University of Aeronautics \& Astronautics\\ \small
Nanjing, Jiangsu, 210016, PR China \\\small $^{4}$ School of Mathematical Sciences, MOE-LSC and SHL-MAC,\\\small
Shanghai Jiao Tong University, 800 Dongchuan road, Shanghai, 200240, P. R. China} \date{} \maketitle
\vspace{3mm}
 {\bf Abstract.} For a  symmetric bivariable   function $f(x,y)$,  let the {\it connectivity function} of a connected graph $G$ be $M_f(G)=\sum_{uv\in E(G)}f(d(u),d(v))$,  where $d(u)$ is the degree of vertex $u$.
 In this paper, we prove that   for an escalating (de-escalating) function $f(x,y)$,  there exists a BFS-graph  with the maximum (minimum) connectivity function $M_f(G)$   among all graphs with a $c-$cyclic degree sequence $\pi=(d_1,d_2, \ldots, d_n)$  and $d_n=1$, and obtain the majorization theorem for connectivity function for unicyclic and bicyclic degree sequences. Moreover, some applications of graph invariants based on degree are included.


  \baselineskip=0.30in
\section{ Introduction}
Throughout this paper,  $G$ denotes  a simple connected
graph with $n$ vertices and $m$ edges, unless
specified otherwise. If  $m=n+c-1$, then $G$ is called a $c$-{\em cyclic
graph}. In particular,  when  $c=0$, $1$ or $2$,  then $G$ is
called a {\em tree}, {\em unicyclic graph} or {\em bicyclic graph}, respectively.  As usual, denote
$N_{G}(v)$  the neighbor set of vertex $v$ in $G$,   let $d_{G}(v)$   be the degree of $v$, and let $N_{G}[v]=N_{G}(v)\cup \{v\}$.  When there is no confusion, we  simplify $N_{G}(v)$, $N_{G}[v]$  and $d_{G}(v)$ as  $N(v)$, $N[v]$ and $d(v)$, respectively. If  $d(v)=1,$ then  $v$ is called  a {\em pendant vertex}. The  nonnegative integer sequence $\pi=(d_{1},d_{2},...,d_{n})$ is called the {\em
degree sequence} of $G$  if  $d_{i}=d(v_{i})$ holds  for   $1\leq i\leq n$, where $V(G)=\{v_1, v_2, \ldots, v_n\}$. Throughout this paper, we  use $d_{i}$ to  denote the $i$-th largest degree of $G$ and we suppose that $d(v_{i})=d_{i}$, where $1\leq i\leq n$.  Let $\Gamma(\pi)$  be
the class of connected   graphs with  degree sequence $\pi$. \par

Among those  vertex-degree-based graph invariants, the {\em Randi\'c index}  $R(G)$ \cite{Ran1} and {\em second Zagreb index} $M_2(G)$ \cite{Gu1972} are two famous topological indices, where $$R(G)=\sum_{uv\in E(G)}\left(d(u)d(v)\right)^{-\frac{1}{2}}, \quad   and \quad  M_{2}(G)=\sum_{uv\in E(G)}\left(d(u)d(v)\right).$$
  Bollob\'as and Erd\H{o}s \cite{De2003} generalized the   concepts of  Randi\'c index and second Zagreb index to the {\em general Randi\'c index}  $W_{\alpha}(G)$, where
$$W_{\alpha}(G)=\sum_{uv\in E(G)}\left(d(u)d(v)\right)^{\alpha}.$$
Similarly with the   general Randi\'c index  of $G$, the {\em general sum-connectivity index}   \cite{Zhou1}  $\chi_{\alpha}(G)$  of $G$  is constructed  as :
$$\chi_{\alpha}(G)=\sum_{uv\in E(G)}\left(d(u)+d(v)\right)^{\alpha},$$
where $\chi_{-\frac{1}{2}}(G)$ is known as the {\em sum-connectivity index} of $G$ \cite{Zhou2},  $\chi_{2}(G)$ is also called  the {\em third Zagreb  index} of $G$ \cite{Sh1}, while $2\chi_{-1}(G)$ is equal  to the {\em harmonic index} of $G$ \cite{Fa1}.
The {\em reformulated Zagreb index} $Z_2(G)$ of $G$ \cite{Mi1} is a slight modification of $\chi_{2}(G)$, where
$$Z_{2}(G)=\sum_{uv\in E(G)}\left(d(u)+d(v)-2\right)^{2}.$$ Another famous  vertex-degree-based graph invariants, called the {\em Atom-Bond connectivity   index} of $G$ \cite{Es1}, is defined  as follows:
$$ABC(G)=\sum_{uv\in E(G)}\sqrt{\frac{d(u)+d(v)-2}{d(u)d(v)}}.$$

For more results on vertex-degree-based graph invariants, one can refer to \cite{Xu1,Gu2013}.  In order to study on  such graph invariants based on adjacency vertex degree,  Wang \cite{Wang1}  recently proposed a new symmetric function.

A symmetric bivariate  function $f(x,y)$  defined on positive  real numbers  is called {\em escalating} (de-escalating) if
\begin{align}\label{21e}f(x_{1},x_{2})+f(y_1,y_2)\geq \ ( resp. \leq)\  f(x_2,y_1) +f(x_1,y_2)\end{align} holds  for any $x_1\geq y_1>0$  and $x_2\geq y_2>0,$
and the inequality  in (\ref{21e}) is strict  if  $x_1>y_1$ and $x_2>y_2$.
  Furthermore, an escalating (de-escalating)  function $f(x,y)$ is called a {\em good escalating  $($de-escalating$)$ function,} if $f(x,y)$ satisfies
  $$\frac{\partial f(x,y)}{\partial x}>0, \ \ \frac{\partial^{2} f(x,y)}{\partial x^{2}}\geq 0,$$
  and
  $$f(x_1+1,x_2)+f(x_1+1,y_1-1)\geq f(x_2,y_1)+f(x_1,y_1)$$   holds for any $x_1\geq y_1$ and $x_2\geq 1$.

Further,  Wang \cite{Wang1} defined the {\em connectivity function} of a connected graph $G$ associated with a symmetric bivariate function   $f(x,y)$  to be
  \begin{align}\label{23e}M_{f}(G)=\sum_{uv\in E(G)}f\left(d(u),d(v)\right),\end{align}
For a connected graph $G$, let us use the notation  $h_G(v)$  to denote  the   distance between  $v$ and  $v_1$, and let $\mathcal {A}_i(G)=\{v:   h_G(v)=i$ and $v\in V(G)\}$. Then,  $\mathcal {A}_0(G)=\{v_1\}$.  If there is no confusion, we simply   write $h(v)$ and  $\mathcal {A}_i$ in place of $h_G(v)$ and  $\mathcal {A}_i(G)$, respectively.
   \begin{definition}\label{21f}
  Let   $G$ be  a connected graph. If  there  exists an ordering $v_{1}\prec v_{2}\prec \cdots\prec v_{n}$ of
  $V(G)$ satisfying the following  $(i)$ and $(ii)$, then $G$ is called a {\em $BFS$-graph} {\em (see \cite{Zhang1})}: \par
 $(i)$    $d(v_{1})\geq d(v_{2})\geq \cdots \geq d(v_{n})$  and $h(v_{1})\leq h(v_{2})\leq\cdots\leq
h(v_{n})$.
 \par  $(ii)$    let $v\in N(u)\backslash N(w)$, $z\in N(w)\backslash N(u)$ such that $h(u)=h(w)=h(v)-1=h(z)-1$, if $u\prec w$, then $v\prec z.$
\end{definition}
In some literatures,   a  $BFS$-tree is also called a {\em greedy tree} (see \cite{Wang1,Ba1}). It is well-known that the $BFS$-tree is unique for any given tree degree sequence $\pi$ \cite{Zhang1}.

In \cite{Wang1},  Wang   proved the following   general important  result.

 \begin{theorem}\label{new0} {\em\cite{Wang1}}
 For an escalating function $f(x,y)$, the connectivity function $W_f(T)$  is maximized by the greedy tree   among trees with given degree sequence.
  \end{theorem}
  Later, Zhang  et al. \cite{Ba1}  studied relations between the extremal trees of two different degree sequences.  In order to state their results, we need the following notation.
Let
$\pi=(x_{1},x_{2},...,x_{n})$ and
$\pi'=(x'_{1},x'_{2},...,x'_{n})$ be  two different
non-increasing  sequences of nonnegative real numbers, we write $\pi\lhd \pi'$
if and only if $\pi\neq \pi'$, $\sum_{i=1}^{n}x_{i}=\sum_{i=1}^{n}x'_{i}$, and
$\sum_{i=1}^{j}x_{i}\leq\sum_{i=1}^{j}x'_{i}$ for all
$j=1,2,...,n$. Such an ordering is sometimes called {\it
majorization} ( for example, see \cite{MA1976}).  Zhang  et al. in \cite{Ba1} proved the following majorization theorem:

\begin{theorem} {\em\cite{Ba1}} Given two tree degree sequences $\pi$ and $\pi^{\prime}$ and an escalating function  with
$\frac{\partial f(x,y)}{\partial x}>0,$ and $\frac{\partial^{2} f(x,y)}{\partial x^{2}}\geq 0,$ let $T_{\pi}^*$ and $T_{\pi^{\prime}}^*$ be the two trees with the maximum connectivity function $M_f(T)$ with degree sequences $\pi$ and $\pi^{\prime}$, respectively.
If  $\pi \lhd \pi^{\prime}$, then  $M_f(T_{\pi}^*)\le M_f(T_{\pi^{\prime}}^*)$.
\end{theorem}

Motivated by the above results, we continued to study on the properties of connectivity function for (good)  escalating or de-escalating functions for $c-$cyclic graphs.


The rest of this paper is organized as follows. In Section 2, we proved that for an escalating (de-escalating) function $f(x,y)$,  there exists a BFS-graph $G$ with the maximum (minimum) connectivity function $M_f(G)$  among all graphs with a $c-$cyclic degree sequence $\pi=(d_1,d_2, \ldots, d_n)$  and $d_n=1$.
In Section 3,  the majorization theorem for connectivity function is obtained for a good escalating function for tree, unicyclic and bicyclic degree sequences. In Section 4, the properties of some graph invariants based on degree are obtained.

\section{The extremal graphs of $\Gamma(\pi)$}
Hereafter, we call $G$ an extremal graph of some graph category  $\mathcal {G}$  if either $M_f(G)$ is maximized in $\mathcal {G}$ when $f(x,y)$ is   escalating   or $M_f(G)$ is minimized in $\mathcal {G}$ when $f(x,y)$ is  de-escalating.

 Let $G-uv$ (respectively, $G-u$)  denote the graph obtained from $G$ by deleting  the   edge $uv\in E(G)$ (respectively, vertex $u$   and the edges incident with it). Similarly, $G+uv$ is a
graph obtained  from $G$ by adding an edge $uv\not\in E(G)$.
\begin{lemma}\label{21l}  Let  $G=(V,E)$  be a  connected  graph with  $u_{1}w_{1}\in
E$, $u_{2}w_{2}\in E$, $w_{1}w_{2}\not\in E$ and
$u_{1}u_{2}\not\in E$, $d(w_{1})\geq d(u_{2})$ and $d(w_{2})\geq d(u_{1})$.   Let
$G'=G-u_{1}w_{1}-u_{2}w_{2}+w_{1}w_{2}+u_{1}u_{2}$. \par
$(i)$ If $f(x,y)$ is escalating, then
$M_{f}(G')\geq M_{f}(G)$, where  $M_{f}(G')> M_{f}(G)$ if and only
if $d(w_{1})>d(u_{2})$ and $d(w_{2})>d(u_{1})$. \par
$(ii)$ If $f(x,y)$ is de-escalating,   then
$M_{f}(G')\leq M_{f}(G)$, where  $M_{f}(G')< M_{f}(G)$ if and only
if $d(w_{1})>d(u_{2})$ and $d(w_{2})>d(u_{1})$.
\end{lemma}
{\bf{Proof.}}    By  (\ref{23e}) and $f(d(u_{1}),d(u_{2}))=f(d(u_{2}),d(u_{1}))$, we have
\begin{align*}
M_{f}(G')-M_{f}(G)=f(d(w_{1}),d(w_{2}))+f(d(u_{2}),d(u_{1}))-f(d(u_{2}),d(w_{2}))-f(d(u_{1}),d(w_{1})).
\end{align*}
  Then  the results follow  from the   definitions of escalating and de-escalating functions. \qed \par \bigskip

Hereafter, let $P_{uv}$ be  a shortest  path  connecting  $u$ and $v$ in $G$, and let $C_q$ be  a cycle with $q$ vertices.
\begin{lemma}\label{22l}
Suppose $G\in \Gamma(\pi)$, and  there exist three vertices $u$, $v$, $w$ of   $G$ such that $uv\in E(G)$, $uw\not\in E(G)$,
$d(v)<d(w)\leq d(u),$ and $d(u)> d(x)$ for all $x\in N(w)$. Then,
$G$ is not an extremal graph of $\Gamma(\pi)$.
\end{lemma}
{\bf{Proof.}}
We firstly suppose that
$uv\not\in P_{uw}$, and we may suppose that $z\in N(w)\cap V(P_{uw})$. If $vz\not\in E(G)$, then let $G_{1}=G+vz+uw-wz-uv$. Clearly, $G_{1} \in \Gamma(\pi)$. Since
$d(u)>d(z)$ and $d(w)>d(v)$,  the result  follows  from Lemma \ref{21l}. Otherwise, $vz\in E(G)$. In this case,
since $d(w)>d(v)$, there exists a vertex $w'\in N(w)$ such that $w'\not\in N[v]$ and $w'\not \in V(P_{uw})$.  Let $G_{2}=G+vw'+uw-ww'-uv$. Then, $G_{2} \in \Gamma(\pi)$. Note that
$d(u)>d(w')$ and $d(w)>d(v)$. Thus,  the result   follows  from   Lemma \ref{21l}.

 We secondly consider the case of  $uv\in P_{uw}$. In this case, since $d(w)>d(v)$, there exists a vertex $w''\in N(w)\setminus V(P_{uw})$ such that  $w''\not\in N[v]$.  Let
$G_{3}=G+vw''+uw-ww''-uv$. Clearly, $G_{3} \in \Gamma(\pi)$ and hence the result   follows  from   Lemma \ref{21l}.   \qed

\begin{lemma}\label{23l} If  $\pi$ is a $c$-cyclic  degree sequence   with $c\geq 0$, then
there exists an extremal graph  $G\in \Gamma(\pi)$   such that  $\{v_{2},v_{3}\}\subseteq N(v_{1})$.
\end{lemma}
{\bf{Proof.}} Let $G$ be an extremal graph of $\Gamma(\pi)$.   If $v_{1}v_{2}\not\in E(G),$  then  there is some vertex $v$ such that $v_{1}v\in E(G)$ and $d(v_{1})\geq d(v_{2})>d(v)$ and $d(v_{1})>d(x)$ holds for all $x\in N(v_{2})$, which contradicts Lemma \ref{22l}. Thus, $v_{1}v_{2}\in E(G).$ Now, we assume that  $v_{1}v_{3}\not\in E(G).$ Then, $d(v_{3})> d(v)$ holds for every $v\in N(v_{1})\setminus\{v_{2}\}$. By  Lemma \ref{22l}, we may assume that there exists some vertex $u\in N(v_{3})$ such that $d(u)=d_{1}$. If $u=v_{2}$, the result already holds. Otherwise, $u\neq v_{2}$, and hence   $d(u)=d_{1}=d_{2}=d_{3}$.   \par

If  $v_{2}\not\in V(P_{v_1 v_3})$,  choose   $u_0\in N(v_{1})\cap V(P_{v_1 v_3})$, since $v_{1}\in N(u_0)\setminus N(v_{3})$, there must  exist  some vertex $v_0\in N(v_{3})\setminus V(P_{v_1 v_3})$ such that $v_0\not\in N[u_0]$. Let $G_{1}=G+v_{1}v_{3}+u_0v_0-v_{1}u_0-v_{3}v_0$.  By Lemma \ref{21l} and  $G_{1}\in \Gamma(\pi)$,   the result already  holds. Otherwise,  $v_{2}\in V(P_{v_1 v_3})$. In this case, $v_1\not\in P_{v_2v_3}$ and $d_1=d_2=d_3$.  Now, it can be proved similarly with  the case  $v_{2}\not\in V(P_{v_1v_3})$. \qed

\par\bigskip
 Let $N_{G}(v, p)$ be the neighbor set of vertices of  $v$ in  $G$ with degree at least  $p$.
\begin{lemma}\label{24l} Let  $C_q$ be a cycle of an extremal graph   $G$ of $\Gamma(\pi)$  with  $w_1w_2\in E(C_q)$ and  $P=u_1\cdots u_{s-1}u_{s}$   being  a  path connecting  $u_1$ and $u_{s}$ such that $u_1\in V(C_q)$ and  $d(w_1)\geq d(w_2)>d(u_s)$. If $d(u_{k})\geq d(w_2)$ holds for some vertex $u_{k}\in V(P)\backslash \{u_1\}$, $N(w_1)\cap \{u_{k},u_{k+1},...,u_s\}=\emptyset$ and $N(w_2)\cap \{u_{k+1},u_{k+2},...,u_s\}=\emptyset$, then there exists an extremal graph $G'$ of $\Gamma(\pi)$ such that  $P'=u_1u_2\cdots u_{k} $ is a part of one cycle in $G'$ with  $N_{G}(u_{k}, d(w_2))\subseteq N_{G'}(u_{k}, d(w_2))$ and $N_{G}(u_{i})= N_{G'}(u_{i})$, where $2\leq i\leq k-1$. Furthermore,   every cycle  of $G$, which do not contain the edge $w_1w_2$, is also a cycle of $G'$.
 \end{lemma}
  {\bf{Proof.}} By the hypothesis, we may suppose that $u_k$ is the last vertex of $P$ such that $d(u_{k})\geq d(w_2)$, namely,  $\max\{d(u_{j}):k+1\leq j\leq s\}<d(w_2).$ Since $d(u_s)<d(w_2)$, we have $2\leq k\leq s-1$. If $P'=u_1u_2\cdots u_k $ is  a part of some   cycle of $G$, then the result already holds by setting  $G'=G$. Thus,   we may suppose that $P'=u_1u_2\cdots u_k$ is not  a part of any  cycle of $G$.  \par

  Let $G_{1}=G+w_1u_{k}+w_2u_{k+1}-w_1w_2-u_{k}u_{k+1}$. Since  $d(w_2)\leq d(u_{k})$ and $d(w_1)\geq d(w_2)> d(u_{k+1})$,  $G_1$ is also an extremal graph of $\Gamma(\pi)$ by Lemma \ref{21l}.  In this case,    $P'=u_1u_2\cdots u_{k}$ is a part of one cycle  of $G_{1}$ such that  $N_{G}(u_{k},d(w_2))\subseteq N_{G_1}(u_{k},d(w_2))$ and   $N_{G}(u_{i})= N_{G_1}(u_{i})$, where $2\leq i\leq k-1$. Furthermore,   every cycle  of $G$, which do not contain the edge $w_1w_2$, is also a cycle of $G_1$. Thus,  the result holds.   \qed

\begin{remark}\label{21r} {\em  By an observation to  the proof of Lemma \ref{24l}, $u_1\in \{w_1, w_2\}$ is also   permitted and  $u_1\not\in \{w_1, w_2\}$ guarantees the existence of   an extremal graph $G'$   of $\Gamma(\pi)$   such that $P'=u_1u_2\cdots u_{k} $ is a part of one cycle in $G'$ with  $N_{G}(u_{k}, d(w_2))\subseteq N_{G'}(u_{k}, d(w_2))$ and $N_{G}(u_{i})= N_{G'}(u_{i})$, where $1\leq i\leq k-1$. }
 \end{remark}

Let $G$ be a connected  graph with $X\subseteq V(G)$. Denote by $G[X]$ the subgraph induced by $X$, and denote by  $\mathcal {R}(G)$ the {\em base graph} obtained  from $G$ by recursively  deleting pendant vertices of the resultant graph until no pendant vertices remain.  If $G$ is a   $c$-cyclic graph, it is easy to see that   $\mathcal {R}(G)$ is  also a  $c$-cyclic graph for    $c\geq 1$   and $\mathcal {R}(G)=\emptyset$ holds for $c=0$.
\begin{lemma}\label{25l} Let  $\pi$ be a  $c$-cyclic degree sequence. If $c\geq 1$ and $d_n=1$,  then there exists an extremal graph   $G$ of $\Gamma(\pi)$  such that     $\{v_{1},v_{2},v_{3}\}$  forms a triangle.
\end{lemma}
{\bf{Proof.}} It suffices to  prove the following  four claims:

{\bf Claim  1.}  There exists an extremal graph   $G$ of $\Gamma(\pi)$  such that   $\{v_{2},v_{3}\}\subseteq N(v_{1})$ and $v_1$ lies on a cycle of $G$.

 By Lemma \ref{23l}, there is  an extremal graph $G\in \Gamma(\pi)$ such that  $\{v_{2},v_{3}\}\subseteq N(v_{1})$. We assume that  claim 1 does not hold.

 If $v_1\not\in V(\mathcal {R}(G))$, then there must exist  one edge $w_1w_2$ in a cycle and one vertex $u$ in $N(v_1)$ such that $u\not\in V(\mathcal {R}(G))$ and $d(u)\geq d(w_1)\geq d(w_2)$ (since $\{v_{2},v_{3}\}\subseteq N(v_{1})$). By Lemma \ref{24l} and Remark \ref{21r}, there exists an extremal graph   $G_1$ of $\Gamma(\pi)$  such that   $\{v_{2},v_{3}\}\subseteq N_{G_1}(v_{1})$ and $v_1$ lies on a cycle of $G_1$,   a contradiction.

 Otherwise, $v_1\in V(\mathcal {R}(G))$. Now, since $v_1$ does not lie on any cycle of $G$, $v_1$ is a cut vertex and  $G-v_1$ contains $d_1$ components, say    $D_1$, $D_2$, ..., $D_{d_1}$. Since $d_n=1$, we may suppose that $v_n$ lies on  $D_1$ and $D_2$ is a component containing a cycle (Recall that $v_1\in  V(\mathcal {R}(G))$, $G-v_1$ contains at least two components containing cycles). In this case, by Lemma \ref{24l} and Remark \ref{21r}, there exists an extremal graph   $G_2$ of $\Gamma(\pi)$  such that   $\{v_{2},v_{3}\}\subseteq N_{G_2}(v_{1})$ and $v_1$ lies on a cycle of $G_2$, a contradiction. This completes the proof of Claim 1.

 {\bf Claim  2.}  There exists an extremal graph   $G$ of $\Gamma(\pi)$  such that   $v_{1}v_{2}$   lies on a cycle of $G$ and $\{v_{2},v_{3}\}\subseteq N(v_{1})$.

 By Claim 1, we may directly  suppose that $G$ is  an extremal graph of  $\Gamma(\pi)$ such that  $\{v_{2},v_{3}\}\subseteq N(v_{1})$ and $v_1$ lies on some cycle $C_q$ of $G$.  We assume that  claim 2 does not hold. Let $d(w_0)=\max\{d(v):v\in V(C_q)\setminus \{v_1\}\}$. If $d(w_0)\geq d(v)$ for some vertex $v\in N(v_2)\backslash \{v_1\}$, then $vw\not\in E(G)$ and $v_2w_0\not\in E(G)$, where $w\in (N(w_0)\cap V(C_q))\setminus \{v_1\}$ (Otherwise, Claim 2 already holds).  Let $G_3=G+v_2w_0+vw-w_0w-v_2v$. Since $d(w_0)\geq d(v)$ and $d(v_2)\geq d(w)$, $G_3$ is an extremal graph of $\Gamma(\pi)$ such that Claim 2 holds, a contradiction. Otherwise, $2\leq d(w_0)<d(v)$ holds for each $v\in N(v_2)$. Let $v'\in N(v)\setminus \{v_2\}$ and let $w'\in (N(v_1)\cap V(C_q))\setminus \{v_3\}$, where $v\in N(v_2)\backslash \{v_1\}$. Since  Claim 2 does not hold,   $v_1v\not\in E(G)$ and $v'w'\not\in E(G).$ Let $G_4=G+v_1v+v'w'-v_1w'-vv'$. Since $d(v_1)\geq  d(v')$ and $d(v)>d(w_0)\geq d(w')$, $G_4$ is an extremal graph of $\Gamma(\pi)$ such that Claim 2 holds. This completes the proof of Claim 2.

{\bf Claim  3.}  There is an extremal graph   $G$ of $\Gamma(\pi)$  such that   $v_{1}v_{2}$ and $v_1v_{3}$ lie on a cycle of $G$.

By Claim 2, we may directly  suppose that $G$ is  an extremal graph of  $\Gamma(\pi)$  such that $v_{1}v_{2}$   lies on a cycle of $G$ and $\{v_{2},v_{3}\}\subseteq N(v_{1})$. Choose   $u_0\in (V(C_{q})\cap N(v_{2}))\setminus\{v_{1}\}$.

 If Claim 3 does not holds, then $v_2v_{3}\not\in E(G)$, $u_0v_3\not\in E(G)$ and there  exists  $v_0\in N(v_{3})\backslash\{v_{1}\}$ such that $v_0\not\in N[u_0]$ (as $d(v_3)\geq d(u_0)\geq 2$), then let $G_{5}=G+v_{2}v_{3}+u_0v_0-v_0v_{3}-u_0v_{2}$. By Lemma \ref{21l}, $G_{5}$ is also an extremal graph of   $\Gamma(\pi)$ such that $\{v_{1},v_{2},v_{3}\}$ forms a triangle of $G_5$, a contradiction.  Thus, Claim 3 holds.

{\bf Claim  4.}  There is an extremal graph   $G$ of $\Gamma(\pi)$  such that   $\{v_{1},v_{2},v_3\}$ forms a triangle  of $G$.

  By Claim 3,  we may directly  suppose that $G$ is an extremal graph of $\Gamma(\pi)$ such that  $\{v_{1}v_{2},v_{1}v_{3}\}$ are two edges of some cycle of $G$, while   $v_{2}v_{3}\not\in E(G)$. Furthermore, we may suppose that
   $$q=\min\big\{p,\hspace{5pt}\text{where}\hspace{5pt}v_{1}v_{2}\in E(C_p)\hspace{5pt}\text{and}\hspace{5pt}v_{1}v_{3}\in E(C_{p})\big\},$$
   and hence $|N(v_2)\cap V(C_q)|=|N(v_3)\cap V(C_q)|=2.$
   We assume that Claim 4 does not hold.   Two cases occur as follows:

  {\bf Case 1. $d_{2}\geq 3$.}

  Choose $w_1\in (N(v_{3})\cap V(C_{q}))\backslash\{v_{1}\}$. If there exists some vertex  $w_2\in N(v_{2})\backslash V(C_{q})$ such that $w_1w_2\not\in E(G)$, then let $G_{6}=G+v_{2}v_{3}+w_1w_2-w_2v_{2}-w_1v_{3}$. By Lemma \ref{21l}, $G_{6}$ is also an extremal graph of   $\Gamma(\pi)$ such that  $\{v_{1},v_{2},v_{3}\}$ forms a triangle in $G_6$, and hence Claim 4 holds, a contradiction.   Otherwise, $w\in N(w_1)$ holds for each   $w\in N(v_{2})\backslash V(C_{q})$, which implies that $d(v_2)=d(v_3)=d(w_1)\geq 3$.

   Nota that $v_1w_1\not\in E(G)$. Otherwise, $\{v_{1},v_{3},w_1\}$ forms a triangle of $G$, and hence Claim 4 holds (since $d(v_2)=d(v_3)=d(w_1)$), a contradiction.  If there exists some vertex $z\in  N(v_{1})\backslash  \{v_2,v_3\}$ such that $zw\not\in E(G)$, then let $G_{7}=G+v_{1}w_{1}+wz-w_1w-v_1z$. By Lemma \ref{21l}, $G_{7}$ is also an extremal graph of   $\Gamma(\pi)$ such that  $\{v_{1},v_{3},w_1\}$ forms a triangle in $G_7$,   a contradiction. Thus,  $z\in N(w)$ holds for each $z\in  N(v_{1})\backslash \{v_2,v_3\}$. Since   $w_1\in N(w)\setminus N(v_1)$, we have  $d(v_1)=d(v_2)=d(v_3)=d(w_1)=d(w)\geq 3$.

     In this case, $wv_3\not\in E(G)$ (Otherwise, Claim 4 already holds). Since $w\in N(v_2)\setminus N(v_3)$ and $d(v_2)=d(v_3)$, there exists some vertex $z'\in N(v_3)\backslash V(C_{q})$ such that $v_2z'\not\in E(G)$. Let $G_8=G+v_{3}w+v_2z'-v_2w-v_3z'$. Since $d(v_1)=d(v_2)=d(v_3)=d(w_1)=d(w)$, by Lemma \ref{21l}, $G_8$ is also an extremal graph of   $\Gamma(\pi)$ such that  $\{v_{3},w_{1},w\}$ forms a triangle of $G_8$, contrary with $d(v_1)=d(v_2)=d(v_3)=d(w_1)=d(w)$.

 {\bf Case 2. $d_2=2$.}

  Suppose $P=v_1u_1\cdots u_s$ is a longest path  among these paths connecting  $v_1$ and a pendant  vertex in $G$. Let $G_9=G+v_{2}v_{3}+u_sv_0-v_0v_{3}-u_0v_{2}$, where  $v_0\in (N(v_{3})\cap V(C_{q}))\setminus \{v_1\}$ and $u_0\in (N(v_{2})\cap V(C_{q}))\setminus \{v_1\}$ (Here, $u_0=v_0$ is also permitted).

If $s\geq 2$, then $M_f(G_9)=M_f(G)$ and hence $G_9$ is also an extremal graph of   $\Gamma(\pi)$ such that $\{v_{1},v_{2},v_{3}\}$ forms a triangle in $G_9$, and hence Claim 4 holds. Otherwise, $s=1$. In this case, $M_f(G_9)-M_f(G)=f(d(v_1),2)+f(2,1)-f(d(v_1),1)-f(2,2)$. By the definition of $f(x,y)$, $G_9$ is also an extremal graph of   $\Gamma(\pi)$ such that $\{v_{1},v_{2},v_{3}\}$ forms a triangle in $G_9$.  \qed

\begin{lemma}\label{26'l} Let $G$ be an extremal $c$-cyclic graph of $\Gamma(\pi)$ with $d_{n}=1$ and $c\geq 1$. If  $\{v_{2},v_{3},...,v_{k}\}\subseteq N_{\mathcal {R}(G)}(v_1)$ and $G[\{v_{1},v_{2},v_{3},...,v_{k}\}]\subseteq  \mathcal {R}(G)$, then   there exists an extremal graph   $G'\in \Gamma(\pi)$  such that   $G'$ is a $BFS$-graph with $\{v_{2},v_{3},...,v_{k}\}\subseteq N_{\mathcal {R}(G')}(v_1)$ and $G[\{v_{1},v_{2},v_{3},...,v_{k}\}]\subseteq  \mathcal {R}(G')$.
\end{lemma}
 \noindent{\bf{Proof.}} We  create an
ordering $\prec$ of $V(G)$ by the  breadth-first-search method as follows: Firstly,  let  $v_{1}\prec v_{2}\prec \cdots\prec  v_{k}$;  secondly,  append all
neighbors   $u_{k+1}$, $u_{k+2}$, ..., $u_{d_{1}+1}$ of $N(v_{1})\setminus\{v_{2},v_{3},...,v_{k}\}$ to the ordered
list, these neighbors are ordered such that $u\prec v$ whenever
$d(u)>d(v)$  (in the remaining case the ordering can be
arbitrary); thirdly  append all
neighbors   of $N(v_{2})\setminus N[v_{1}]$ to the ordered
list, these neighbors are ordered such that $u\prec v$ whenever
$d(u)>d(v)$  (in the remaining case the ordering can be
arbitrary). Then, with the same method we can append  the vertices
$N(v_{i})\setminus (N[v_{1}]\cup N(v_{2})\cup \cdots\cup N(v_{i-1}))$ in the ordered list, where $3\leq i\leq k$, and then to  the
vertices $N(w)\setminus(N[v_{1}]\cup N(v_{2})\cup \cdots\cup N(v_{k}))$, where   $d(w)=\max\big\{d(z):$    $z\in N(v_{1})\setminus\{v_{2},v_{3},...,v_k\}\big\}$.
 Continue recursively with
all vertices $v_{1},v_{2},...,  v_n$ until all vertices of $G$ are
processed.    It suffices to show the following two   claims. \par
{\bf Claim 1.} If $i<j$, then  $d(u)\geq d(v)$ holds for any $u\in \mathcal {A}_i$ and $v\in \mathcal {A}_j$.
\par
Suppose that there exist vertices  $u$ and $v$ with  $h(u)<h(v)$, while   $d(u)<d(v)$.    Furthermore, for convenience, we may choose  $u$ as the first vertex in the ordering $\prec$ with such property and let $d(v)=\max\big\{d(w):$  where  $h(w)>h(u)\big\}$.    Then, $1\leq i<j$  and  $u\not\in \{v_1,v_2,...,v_k\}$. \par

If  $u\in V(P_{vv_{1}})$, then we choose $u'\in \mathcal {A}_{i-1}\cap V(P_{vv_{1}})$ such that $u'u\in E(G).$ Since $d(v)>d(u)$, there exits vertex $v'\in N(v)\setminus V(P_{vv_{1}})$ such that $v'\not\in N(u)$.  By the choice of $u$ and $v'$, we have  $d(u')\geq d(v')$ and $d(v)>d(u)$. Let $G_{1}=G+uv'+u'v-u'u-vv'.$ By  Lemma \ref{21l} and  the choice of $G$, we can conclude that $G_1$ is also an extremal graph of  $\Gamma(\pi)$.
Now, we construct a new ordering $\prec'$ of $V(G_{1})$ with the similar method as $\prec$. We suppose that $v_{1}\prec v_{2}\prec \cdots\prec v_{k}\prec u_{1}\prec u_{2}\prec \cdots \prec u_{t}\prec u$ are the first $t+k+1$ elements   in the ordering $\prec$ of $V(G)$. By the choice of $u$ and $v$, we have $v_{1}\prec' v_{2}\prec' \cdots\prec' v_{k}\prec' u_{1}\prec' u_{2}\prec' \cdots \prec' u_{t}\prec' v$ are the first $t+k+1$ elements  in the ordering $\prec'$ of $V(G_{1})$. By the choice of $u$ and $v$,   for every $w\in \{v_{1}, v_{2},...,v_{k},u_{1},u_2,  ..., u_{t},v\}$, if $w\in \mathcal {A}_{t}(G_{1})$, then $d_{G_{1}}(w)\geq d_{G_{1}}(w')$ holds for every $w'\in \mathcal {A}_{s}(G_{1})$, where $t<s$.

In what follows, we suppose that  $u \not \in V(P_{vv_{1}})$.

 Choose  $w\in V(P_{vv_{1}})\cap \mathcal {A}_{j-1}$ and choose $u'\in \mathcal {A}_{i-1}$ such that $wv\in E(G)$ and  $u'u\in E(G)$. By the former argument, we may suppose that  $u\neq w.$

If   $uw\not\in E(G)$, by the  choice of $u$ and $j-1>i-1$, $d(u')\geq d(w)$. Recall that $d(v)>d(u).$  Let $G_{2}=G+uw+u'v-uu'-vw.$ By Lemma \ref{21l},  $G_2$ is also an extremal graph of  $\Gamma(\pi)$. \par

Otherwise,  $uw\in E(G)$.  Recall that $d(v)>d(u).$  There exists vertex $w'\in N(v)$ such that $w'\not\in V(P')\cup N[u]$, where $P'$ is a path connected $v$ and $u$ with $\{u',v_1,w\}\subseteq V(P')$.  Let $G_{2}=G+uw'+u'v-uu'-vw'.$ By the  choice of $u$ and $j-1>i-1$, $d(u')\geq d(w')$. In this case,   Lemma \ref{21l} implies that  $G_2$ is also an extremal graph of  $\Gamma(\pi)$.

 Now, we construct a new ordering $\prec'$ of $V(G_{2})$ with the similar method as $\prec$. We suppose that $v_{1}\prec v_{2}\prec \cdots\prec v_{k}\prec u_{1}\prec u_{2}\prec \cdots \prec u_{t}\prec u$ are the first $t+k+1$ elements   in the ordering $\prec$ of $V(G)$. By the choice of $u$ and $v$, we have $v_{1}\prec' v_{2}\prec' \cdots\prec' v_{k}\prec' u_{1}\prec' u_{2}\prec' \cdots \prec' u_{t}\prec' v$ are the first $t+k+1$ elements  in the ordering $\prec'$ of $V(G_{2})$. By the choice of $u$ and $v$,   for every $w\in \{v_{1}, v_{2},...,v_{k},u_{1},u_2,  ..., u_{t},v\}$, if $w\in \mathcal {A}_{t}(G_{2})$, then $d_{G_{2}}(w)\geq d_{G_{2}}(w')$ holds for every $w'\in \mathcal {A}_{s}(G_{2})$, where $t<s$.

  Repeating the above process finitely many times, we can achieve a   graph $G^{*}$  and an  ordering $\prec^{*}$ such that Claim 1 holds for $G^{*}$. Thus, we may directly suppose that Claim 1 holds for $G$.
\par

{\bf Claim 2.} If $h(u')=h(v')=i$ and $u'\prec v'$, then  $d(u)\geq d(v)$ holds for every $u\in (N(u')\cap \mathcal {A}_{i+1})\setminus N(v')$ and every $v\in (N(v')\cap  \mathcal {A}_{i+1})\setminus N(u').$ \par

   If Claim 2 does  not hold, we may suppose that $u'$ is the first vertex  in the ordering $\prec$ with the property that there exist vertices $v'\in \mathcal {A}_{i}$ (suppose that  $u'\in \mathcal {A}_{i}$), $u\in (N(u')\cap \mathcal {A}_{i+1})\setminus N(v')$ and  $v\in (N(v')\cap \mathcal {A}_{i+1})\setminus N(u')$   such that $u'\prec v'$, but $d(u)<d(v).$ Furthermore, according to this, we may suppose that $u$ is the first vertex in the ordering $\prec$ and  suppose that $d(v)=\max\{d(z): z\in \mathcal {A}_{i+1}\setminus N(u'),$ where $u\prec z\}$. It is easy to see that  $u\not\in \{v_1,v_2,...,v_k\}$. Let $G_{3}=G+u'v+v'u-u'u-v'v.$ Then, $G_{3}\in \Gamma(\pi)$. By the choice of $u$, we have  $d(u')\geq d(v')$.   Thus, Lemma \ref{21l} implies that  $G_{3}$ is also an extremal graph of $\Gamma(\pi)$. Now, we construct a new ordering $\prec'$ of $V(G_{3})$ with the similar method as $\prec$.
     We suppose that $v_{1}\prec v_{2}\prec \cdots\prec v_{k}\prec u_{1}\prec u_{2}\prec \cdots \prec u_{t}\prec u$ are  the first $t+k+1$   elements   in the ordering $\prec$ of $V(G)$. By the choice of $u$ and $v$,  $v_{1}\prec' v_{2}\prec' \cdots\prec' v_{k}\prec' u_{1}\prec' u_{2}\prec' \cdots \prec' u_{t}\prec' u$ are  the first $t+k+1$ elements  in the ordering $\prec'$ of $V(G_{3})$. By the choice of $u$ and $v$,   for every $z\in \{v_{1}, v_{2},...,v_{k},u_{1},u_2,..., u_{t},v\}\cap \mathcal {A}_{i+1}(G_{3})$, if $w\in \mathcal {A}_{i+1}(G_{3})$ and $z\prec' w$, then $d_{G_{3}}(z)\geq d_{G_{3}}(w)$.

     Repeating the above process finitely many times, we can obtain a   graph $G^{**}$ and an ordering $\prec^{**}$ of $V(G^{**})$ such that $M_{f}(G^{**})=M_{f}(G)$ and  Claim 2 holds for $G^{**}$.  \par  This completes the proof of this result.   \qed

\begin{theorem}\label{21t}  If $d_{n}=1$, then   there exists an extremal graph   $G\in \Gamma(\pi)$  such that   $G$ is a $BFS$-graph. Furthermore,  $\{v_{1},v_{2},v_{3}\}$ forms a triangle of $G$ when $c\geq 1$.
\end{theorem}
 \noindent{\bf{Proof.}} We may suppose that $c\geq 1$, as the case of $c=0$ can be proved easily with the same argument.  By Lemma \ref{25l}, we may suppose that    $G$ is an extremal graph in $\Gamma(\pi)$ such that   $G$ contains a triangle with $V(C_{3})=\{v_{1},v_{2},v_{3}\}$. Now,  the result follows from Lemma \ref{26'l}.   \qed

\begin{figure}[H]
\vspace*{-1.9cm}\begin{center} \includegraphics[scale=0.45]{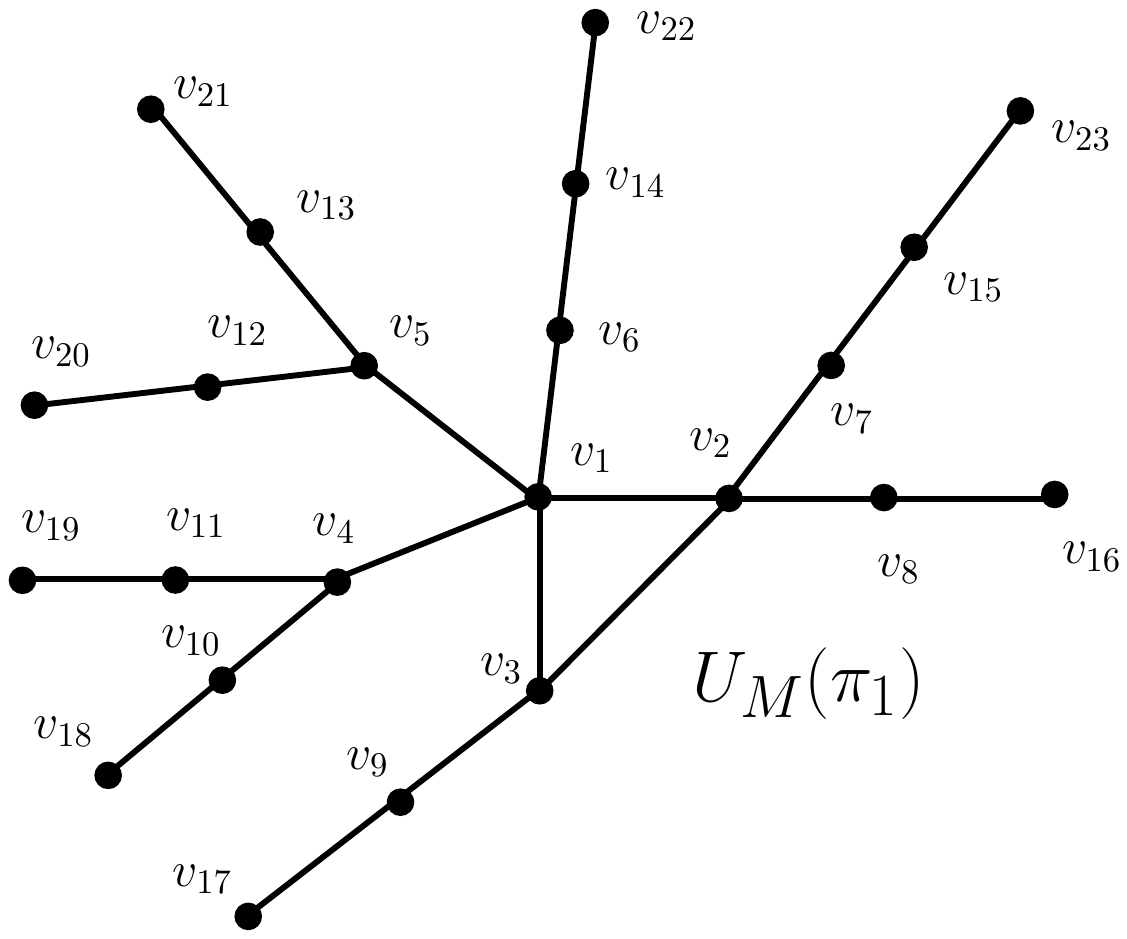} \par
\vspace*{-6.7cm}
\caption{The unicyclic   graph    $U_{M}(\pi_{1})$.} \label{liu1f}\end{center}
     \end{figure} \par \vspace{-0.8cm}
 Hereafter, we use the symbol   $p^{(q)}$ to define  $q$   copies  of the real  number $p$.
Suppose $\pi=(d_{1},d_{2},...,d_{n})$, where $d_{n}=1$. It has been shown that the $BFS$-tree is unique for any tree degree sequence $\pi$ \cite{Zhang1}. Actually, we can construct a unique unicyclic $BFS$-graph   $U_{M}(\pi)$ by the following  breadth-first-search method for any unicyclic  degree sequence $\pi$:  The unique cycle  of  $U_{M}(\pi)$ is a triangle with $V(C_{3})=\{v_{1},v_{2},v_{3}\}$. Select the  vertex $v_{1}$  as the  root vertex and begin with $v_{1}$ of the zeroth layer. Select the vertices $v_{2}$, $v_{3}$, $v_{4}$, $v_{5}$, ..., $v_{d_{1}+1}$ as the  first  layer such that $N(v_{1})=\{v_{2},v_{3},v_{4},v_{5},...,v_{d_{1}+1}\}$. Let   $N(v_{2})=\{v_{1},v_{3},v_{d_{1}+2},v_{d_{1}+3},...,v_{d_{1}+d_{2}-1}\}$ and  $N(v_{3})=\{v_{1},v_{2},v_{d_{1}+d_{2}},...,v_{d_{1}+d_{2}+d_{3}-3}\}$. Then, append $d_{4}-1$ vertices to $v_{4}$ such that  $N(v_{4})=\{v_{1},$$v_{d_{1}+d_{2}+d_{3}-2},...,v_{d_{1}+d_{2}+d_{3}+d_{4}-4}\}$ $\cdots$.
Informally, for a given  unicyclic degree sequence $\pi_{1}=(5,4,3^{(3)},2^{(10)},1^{(8)})$, $U_{M}(\pi_{1})$ is the unicyclic graph   as  shown in Figure \ref{liu1f}.

\begin{corollary}\label{21c} Let $\pi$ be a $c$-cyclic degree sequence with $d_{n}=1$.\par $(i)$  {\em\cite{Wang1}} If $c=0$, then  the $BFS$-tree is  an extremal   graph  of  $\Gamma(\pi)$.\par $(ii)$ {\em\cite{Zhangc}} If  $c=1$, then $U_{M}(\pi)$   is  an extremal   graph  of  $\Gamma(\pi)$.
\end{corollary}
\noindent{\bf{Proof.}} One can easily check that  $(i)$  and $(ii)$ follows from Theorem  \ref{21t}.    \qed

\par
\begin{figure}[H]
\vspace*{-1.9cm}\begin{center} \includegraphics[scale=0.65]{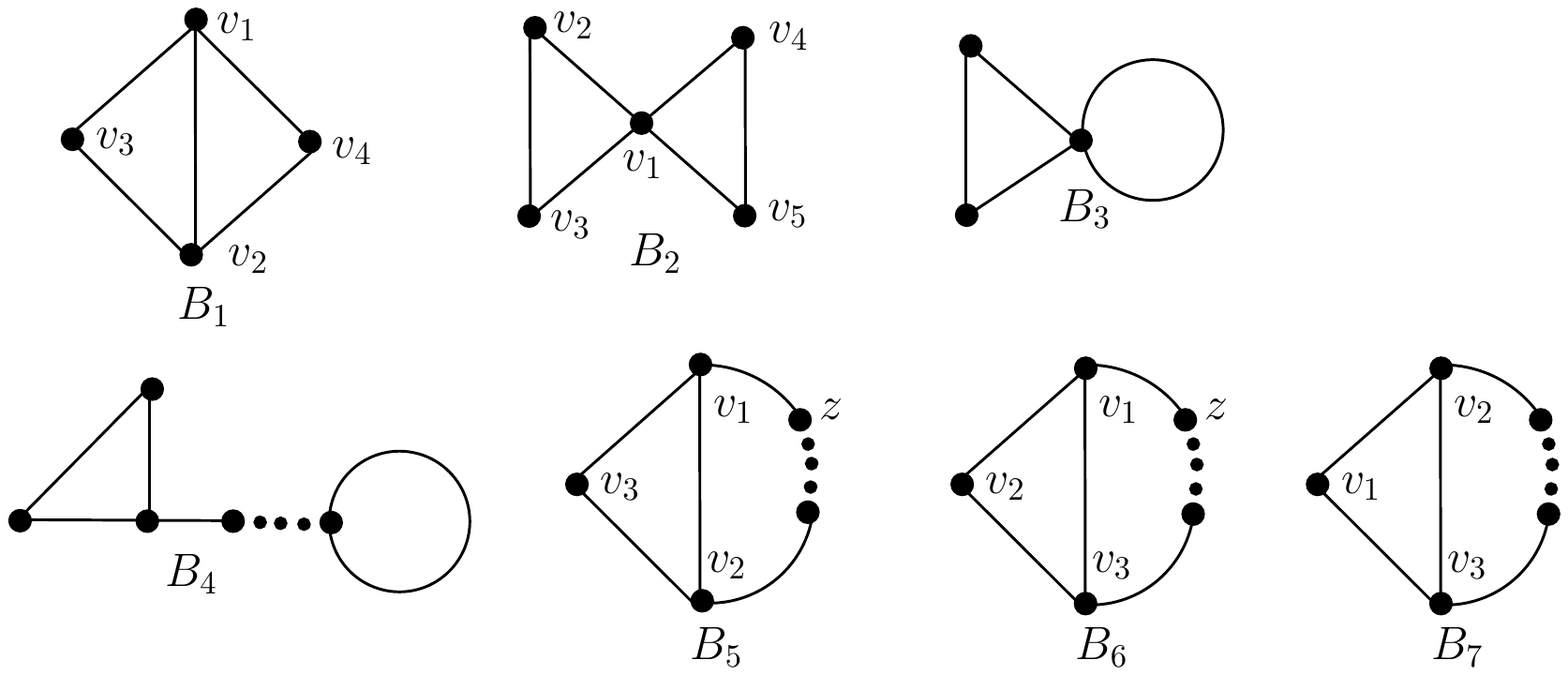} \par\vspace*{-11.5cm}
\caption{The bicyclic   graphs    $B_1$, $B_2$, $...$,  $B_7$.} \label{liu2f}\end{center}
     \end{figure}

\bigskip
Paths $P_{l_{1}}$, $P_{l_{2}}$, $...$, $P_{l_{k}}$
are said to {\it have almost equal lengths} if $l_{1}$, $l_{2}$, $...$,
$l_{k}$ satisfy $|l_{i}-l_{j}|\leq 1$ for $1\leq i\leq j\leq k$.  In what follows, let $B_1$, $B_2$, $...$,   $B_7$ be seven bicyclic graphs as shown in Figure  \ref{liu2f}.  If $\pi=(d_{1},d_{2},...,d_{n})$ is a bicyclic degree sequence, then   $\sum_{i=1}^{n}d_{i}=2n+2$, which implies that    $\pi$ should
be one of the following four cases.  Moreover, we construct a special bicyclic graph $B_{M}(\pi)$ of $\Gamma(\pi)$ as follows:\par
$(i)$ If $d_{n} = 1$ and $d_{1}\geq  d_{2}\geq 3$, let $B_{M}(\pi)$ be a   $BFS$-graph such that $\mathcal {R}(B_{M}(\pi))\cong B_1$ and the
remaining vertices appear in a    $BFS$-ordering.
 \par
 $(ii)$  If $d_{1}\geq 5>d_{2}= 2$ and $d_{n} = 1$, let $B_{M}(\pi)$ be the  bicyclic  graph with  $n$  vertices obtained from $B_2$ by attaching  $d_{1}-4$ paths of almost equal lengths to the maximum degree vertex of $B_2$.
\par
$(iii)$ If $\pi=(4,2^{(n-1)})$,   let   $B_{M}(\pi)\cong
B_3$.\par
$(iv)$ If $\pi=(3^{(2)},2^{(n-2)})$, let  $B_{M}(\pi)\cong
B_5$.\par

 \par \vspace{-0.8cm}
\begin{lemma}\label{27l'} If $\pi$ is a bicyclic degree sequence with  $d_1\geq d_2\geq 3$ and  $d_{n}=1$, then $B_{M}(\pi)$ is an extremal bicyclic graph of $\Gamma(\pi)$.
\end{lemma}
\noindent{\bf{Proof.}}    By Theorem \ref{21t}, we may suppose that $G$ is an extremal bicyclic graph of $\Gamma(\pi)$ such that $\{v_{1},v_{2},v_{3}\}$ forms a triangle $\Delta_1$ of $G$ and $G$ is an $BFS$-graph. We firstly prove  the following two claims.  \par
{\bf Claim 1.} There is an  extremal bicyclic graph $G$ of $\Gamma(\pi)$ such that  $\{v_{1},v_{2},v_{3}\}$ forms a triangle $\Delta_1$  of $G$  and $G$ contains a cycle $C_q$ with $\{v_1v_2,v_1v_4\}\subseteq  E(C_q)$.

 Since $G$ is a bicyclic graph, there is another cycle, say $C_q$ such that  $C_q\neq \Delta_1$.  If $v_1v_2$ does not lie  on   $C_q$, then $\mathcal {R}(G)\cong B_3$ or $\mathcal {R}(G)\cong B_4$. It is easy to see that  there exists at least one edge $uv$ on $C_q$ such that $d(u)\geq d(v)$, $\{v_1,v_2,v_3\}\cap\{u,v\}=\emptyset$ and $d_5\geq d(v)$, by Lemma  \ref{24l} and Remark \ref{21r}, there exits another  extremal graph $G'$ such that $\Delta_1$  is also  a triangles of $G'$,  $v_1v_2$ also lies on   a cycle of $G'$ and $\textcolor{red}{\{v_2,v_3,v_4\}}\subseteq N_{G'}(v_1).$  Thus, we may suppose that $\Delta_1$ is a triangle  of $G$, $v_1v_2$   lies on the   cycle $C_q$ of $G$ and $\textcolor{red}{\{v_2,v_3,v_4\}}\subseteq N_{G}(v_1).$ In this case,   $\mathcal {R}(G)\in \big\{ B_5,\,B_6,\,B_7\big\}$.

If either  $v_4\in V(C_q)$ or $N(v_4)\cap \{v_{1},v_{2},v_{3}\}\neq \{v_1\}$, then  Claim 1 already holds. Thus, we may suppose that  $v_4\not\in V(C_q)$ and $N(v_4)\cap \{v_{1},v_{2},v_{3}\}=\{v_1\}$ (In this case, $v_4\not\in V(\mathcal {R}(G))$). Let $u_0\in N(v_4)\setminus \{v_{1},v_{2},v_{3}\}$ (such vertex exists since $G$ contains a cycle $C_q$).

If $\mathcal {R}(G)\cong B_5$ or $\mathcal {R}(G)\cong B_7$, we choose  $u\in N(v_2)\cap V(B_5)$ such that $u\not \in \{v_1,v_3\}$. Let $G_1=G+uu_0+v_2v_4-v_2u-u_0v_4$. Since $d(v_2)\geq d(u_0)$ and $d(v_4)\geq d(u)$,  by Lemma \ref{21l},   $G_1$ is also an extremal graph of $\Gamma(\pi)$  such that Claim 1 holds.

If $\mathcal {R}(G)\cong B_6$, we choose  $u\in N(v_3)\cap V(B_6)$ such that $u\not \in \{v_1,v_2\}$. Let $G_2=G+uu_0+v_3v_4-v_3u-v_4u_0$. Since $d(v_3)\geq d(u_0)$ and $d(v_4)\geq d(u)$,  by Lemma \ref{21l},   $G_2$ is also an extremal graph of $\Gamma(\pi)$  such that Claim 1 holds.   This completes the proof of Claim 1.

 \par
 {\bf Claim 2.} There is an extremal bicyclic graph $G$ of $\Gamma(\pi)$ such that   $\{v_1,v_2,v_{4}\}$ and $\{v_1,v_2,v_{3}\}$ form two triangles of $G$.

By Claim 1, we may suppose that $G$ is an extremal graph  of $\Gamma(\pi)$ such that $\{v_1,v_2,v_{3}\}$  form a triangle  of $G$ and $G$ contains a cycle $C_q$ with $\{v_1v_2,v_1v_4\}\subseteq  E(C_q)$. In this case, either $\mathcal {R}(G)\cong B_5$ or $\mathcal {R}(G)\cong B_6$ with $z=v_4$.   

If $v_3v_4\in E(G)$, then   $\{v_1,v_2,v_3,v_4\}$ form a $C_4$ of $G$ (In this case, $\mathcal {R}(G)\cong B_6$). In this case,  there exists some vertex $u\in N(v_2)\setminus V(\mathcal {R}(G))$ and $v_2v_4\not\in E(G)$. Let $G_3= G+v_2v_4+v_3u-v_3v_4-v_2u$. By Lemma \ref{21l},  $G_3$ is also an extremal graph of $\Gamma(\pi)$  such that Claim 2 holds.

Otherwise, $v_3v_4\not\in E(G).$

If   $\mathcal {R}(G)\cong B_6$, then there exists a vertex $w_2\in N(v_2)\setminus V(\mathcal {R}(G))$. Let $G_4=G+v_2v_4+w_2w_0-v_2w_2-v_4w_0$, where $w_0\in (N(v_4)\cap V(B_6))\setminus\{v_1\}$. By Lemma \ref{21l},  $G_4$ is   an extremal graph  of $\Gamma(\pi)$ such that   Claim 2 holds.

Therefore, we may suppose that   $\mathcal {R}(G)\cong B_5$ in what follows. Furthermore, we may suppose that $v_2v_4\not\in E(G).$ Otherwise, Claim 1 already holds. We choose $u_0\in (N(v_4)\cap V(B_5))\setminus\{v_1\}$ and $z_0\in (N(v_2)\cap V(B_5))\setminus\{v_1,v_3\}$, where $u_0=z_0$ is permitted.

If $d_3\geq 3$, then we choose  $w\in N(v_3)\backslash V(B_5)$. Note that $wu_0\not\in E(G)$.   Let
$G_5=v_3v_4+wu_0-v_3w-v_4u_0$. By Lemma \ref{21l}, $G_5$ is an extremal graph  of $\Gamma(\pi)$ such that  $\{v_1,v_2,v_3,v_4\}$ form a $C_4$ and $\{v_1,v_2,v_3\}$ form a triangle  of $G_5$. From the above argument, Claim 2 holds  (see the case of $v_3\in N(v_4)$ and $\mathcal {R}(G)\cong B_6$).

Otherwise, $d_3=2$, and hence $d_1\geq 4$ (since $d_n=1$ and $\sum_{i=1}^{n}d_i=2(n+1)$). Let $w_1$ be a vertex of $N(v_1)\backslash \{v_2,v_3,v_{4}\}$, and let $G_6=G+v_2v_4+z_0v_1+w_1u_0-v_2z_0-v_4u_0-v_1w_1$.
Now, by  (\ref{23e}) it follows that $M_f(G_6)-M_f(G)=f(d_1,2)+f(2,d(w_1))-f(d_1,d(w_1))-f(2,2),$ and hence $G_6$ is also an extremal graph  of $\Gamma(\pi)$ such that   Claim 2 holds. \par

Now, by Claim 2 and Lemma \ref{26'l} we have  $G\cong B_{M}(\pi)$, as desired.
\qed
\begin{theorem}\label{22t}   If $\pi$ is a bicyclic degree sequence, then $B_{M}(\pi)$   is  an extremal   graph  of  $\Gamma(\pi)$.
\end{theorem}
\noindent{\bf{Proof.}} Since  $\pi=(d_{1},d_{2},...,d_{n})$ is a bicyclic degree sequence, then   $\sum_{i=1}^{n}d_{i}=2n+2$, and hence it suffices to consider the following four cases. \par
 If $\pi=(4,2^{(n-1)})$ (respectively, $(3^{(2)},2^{(n-2)})$), it is easy to check that $B_3$ (respectively, $B_5$) is an extremal   graph  of  $\Gamma(\pi)$. If $d_1\geq d_2\geq 3$ and  $d_{n}=1$, then the result follows from Lemma \ref{27l'}. \par
 Otherwise,   $d_{1}\geq 5>d_{2}= 2$ and $d_{n} = 1$. In this case, by Theorem \ref{21t}, we may suppose that $G$ is an extremal bicyclic  graph of $\Gamma(\pi)$ such that $G$ is a  $BFS$-graph obtained by attaching   $d_{1}-4$ paths of almost equal lengths to the maximum degree vertex of $B_3$. Now, we assume that $B_3\not\cong B_2$, and hence $B_3$ contains a cycle $C_q$, where $q\geq 4$. We suppose that $N(v_1)\cap V(C_q)=\{u_1,w_1\}$, $N(u_1)\cap V(C_q)=\{v_1,u_2\}$ and $N(w_1)\cap V(C_q)=\{v_1,w_2\}$, where $w_2=u_2$ is permitted. We choose $z$ as a pendant  vertex such that the distance between $v_1$ and $z$ is as small as possible, and suppose that $z_1\in N(z)$.

Let $G_1=G+u_1w_1+zw_2-u_1u_2-w_1w_2$.

If $z_1=v_1$, then  $M_f(G_1)-M_f(G)=f(d(v_1),2)
+f(2,1)-f(d(v_1),1)-f(2,2)$. If $z_1\neq v_1$, then  $M_f(G_1)=M_f(G)$. In both case, $G_1$ is also an extremal bicyclic graph of $\Gamma(\pi)$ such that $\mathcal {R}(G_1)\cong B_2$. Now, the result follows from Lemma \ref{26'l}. \qed

\section{The  majorization theorem of connectivity function}
 If $\pi= (d_{1}, d_{2}, ..., d_{n})$ is a non-increasing  integer sequence and $d_{i}\geq d_{j}+2$ holds for some integers $1\leq i<j\leq n$, then the following
operation is called a {\em unit transformation} from $i$ to $j$ on $\pi$: subtract 1 from $d_{i}$ and add 1
to $d_{j}$.  The following famous  lemma about majorization for integer sequences is due to   Muirhead (for instance, see \cite{MA1976,Mh1}).
\begin{lemma}\label{26l}  $($Muirhead Lemma$)$ If $\pi$ and $\pi'$ are two non-increasing  integer sequences and $\pi\lhd \pi'$, then $\pi$ can be obtained from $\pi'$ by a finite sequence of unit transformations.
\end{lemma}

\begin{lemma}\label{27l}
 Let $u$, $v$ be two vertices of a connected   graph   $G$ with $d(u) \geq d(v)\geq 2$
and   $w$ be some vertex of $ N(v)\setminus (V(P_{uv})\cup N[u])$.   Let
 $G'=G+wu-wv$. Suppose that   $N(v)\backslash\{w,u\}=\{z'_{1},z'_{2},...,z'_{|N(v)\backslash\{w,u\}|}\}$, and  there exist   $\{z_{1},z_{2},...,z_{|N(v)\backslash\{w,u\}|}\}\subseteq N(u)\setminus\{v\}$ such that $d(z_{i})\geq d(z'_{i})$ for $1\leq i\leq |N(v)\backslash\{w,u\}|$.   If  $f(x,y)$ is a good  escalating function, then $M_{f}(G')> M_{f}(G).$
\end{lemma}
{\bf{Proof.}}
 If
$uv\not\in E(G)$, by  (\ref{23e})  it follows that
\begin{align*}
&M_{f}(G')-M_{f}(G)\\=&\sum_{z\in
N(u)}\left(f\left(d(u)+1,d(z)\right)-f\left(d(u),d(z)\right)\right)-\sum_{z'\in
N(v)\setminus\{w\}}\left(f\left(d(v),d(z')\right)-f\left(d(v)-1,d(z')\right)\right)
\end{align*}\begin{align*}\hspace{-200pt}+f\left(d(u)+1,d(w)\right)-f\left(d(v),d(w)\right).
\end{align*}

Let $g(x,y)=\frac{\partial f(x,y)}{\partial x}$.
We may suppose that $N(v)\backslash\{u,w\}=\{z'_{1},z'_{2},...,z'_{d(v)-1}\}$. Then, there exists  $\{z_{1},z_{2},...,z_{d(v)-1}\}\subseteq N(u)\setminus \{v\}$ such that $d(z_{i})\geq d(z'_{i})$. Recall that $d(u)\geq d(v)$, $d(z_{i})\geq d(z'_{i})$  and $\frac{\partial^{2} f(x,y)}{\partial x^{2}}\geq 0$.  Thus, for $1\leq i\leq d(v)-1,$ there exists $0\leq a\leq 1$ and $0\leq b\leq 1$ such that
\begin{align*}&f(d(u)+1,d(z_i))-f(d(u),d(z_i))-(f(d(v),d(z'_i))-f(d(v)-1,d(z'_i)))\\= &f(d(u)+1,d(z_i))-f(d(u),d(z_i))-(f(d(v),d(z_i))-f(d(v)-1,d(z_i)))\\&+f(d(v),d(z_i))-f(d(v)-1,d(z_i))-(f(d(v),d(z'_i))-f(d(v)-1,d(z'_i)))
\\\geq &f(d(u)+1,d(z_i))-f(d(u),d(z_i))-(f(d(v),d(z_i))-f(d(v)-1,d(z_i)))\\=& g(d(u)+a,d(z_i))-g(d(v)-1+b,d(z_i))\geq 0,\end{align*}
 which implies that  $M_{f}(G')-M_{f}(G)\geq f\left(d(u)+1,d(w)\right)-f\left(d(v),d(w)\right)>0$, as $g(x,y)>0$.

Otherwise, $uv\in E(G)$. In this case,
 by  (\ref{23e})  it follows that
\begin{align*}
&\hspace{10pt}M_{f}(G')-M_{f}(G)\\=&\sum_{z\in
N(u)\setminus\{v\}}\left(f\left(d(u)+1,d(z)\right)-f\left(d(u),d(z)\right)\right)-\sum_{z'\in
N(v)\setminus\{u,w\}}\left(f\left(d(v),d(z')\right)-f\left(d(v)-1,d(z')\right)\right)
\\&\hspace*{10pt} +f\left(d(u)+1,d(w)\right)-f\left(d(v),d(w)\right)+f\left(d(u)+1,d(v)-1\right)-f\left(d(u),d(v)\right).\end{align*}
Since $|N(u)\setminus\{v\}|>|N(v)\setminus\{u,w\}|$, there exists some vertex $z_0\in
N(u)\setminus \{v\}$ such that  \begin{align*}&\sum_{z\in
N(u)\setminus\{v\}}\left(f\left(d(u)+1,d(z)\right)-f\left(d(u),d(z)\right)\right)-\sum_{z'\in
N(v)\setminus\{u,w\}}\left(f\left(d(v),d(z')\right)-f\left(d(v)-1,d(z')\right)\right)\\\geq&f\left(d(u)+1,d(z_0)\right)-f\left(d(u),d(z_0)\right)>0.\end{align*}
 Furthermore, since  $f(x,y)$ is good escalating,   we can conclude that
\begin{align*}
M_{f}(G')-M_{f}(G)>&f\left(d(u)+1,d(w)\right)-f\left(d(v),d(w)\right)+f\left(d(u)+1,d(v)-1\right)-f\left(d(u),d(v)\right)
\\\geq&0.\end{align*} This completes the proof of this result.
 \qed
\par\bigskip
Suppose that  $\pi=(d_{1},d_{2},\ldots, d_{n})$ and $\pi'=(d'_{1}, d'_{2}, \ldots, d'_{n})$ are two $c$-cyclic graphic degree  sequences.  If $\pi\lhd \pi'$, by Lemma \ref{26l}, we may suppose that  $\pi$ and $\pi'$ differ only in two positions
where the difference is  $1$, that is, $d_{i}=d'_{i}$,
$i\neq p,q$, $1\leq p<q\leq n$, and $d'_{p}=d_{p}+1$,
$d'_{q}=d_{q}-1$.

Let $G$ and $G'$ be  an extremal  $c$-cyclic graph  of  $\Gamma(\pi)$ and $\Gamma(\pi')$, respectively.  If  $w$ is a vertex of $G$  such that $w\in N(v_{q})\setminus (V(P_{v_pv_q})\cup N[v_{p}])$, and there exists  $\{z_{1},z_{2},...,z_{|N(v_q)\backslash\{w,v_p\}|}\}\subseteq N(v_p)\setminus\{v_q\}$ such that $d(z_{i})\geq d(z'_{i})$ for $1\leq i\leq |N(v_q)\backslash\{w,v_p\}|$, where   $N(v_q)\backslash\{w,v_p\}=\{z'_{1},z'_{2},...,z'_{|N(v_q)\backslash\{w,v_p\}|}\}$,  then we call $w$ a {\em  surprising vertex} of $G$. If $f(x,y)$ is   good escalating and    $G$ contains  some surprising vertex  $w$, then let $G_{1}=G+v_{p}w-v_{q}w$. Thus,  $G_{1}\in
 \Gamma(\pi')$.    Since $p<q$,  $d(v_{p})\geq d(v_{q})$ follows from Theorem \ref{21t}. By Lemma
\ref{27l}, we have
 $M_f(G)< M_f(G_{1})\leq M_f(G')$. Thus,
\begin{equation}\label{320e} \text{if}\,\,G\,\,\text{contains a surprising vertex, then}\,\,M_f(G)<M_f(G').\end{equation}

Zhang et al. \cite{Ba1} showed  that the size of  $M_{f}(T)$ and $ M_{f}(T')$ can be deduced from the relation $\pi\lhd \pi'$ when $f(x,y)$ is a special symmetric bivariable  function  for any  extremal trees $T$ and $T'$  of  $\Gamma(\pi)$ and  $\Gamma(\pi')$, respectively. Actually, if $\pi\lhd \pi'$ and $f(x,y)$ is  good  escalating,  by Corollary  \ref{21c} $(i)$ we can conclude that the $BFS$-tree of $\Gamma(\pi)$ contains a surprising vertex. Thus, by (\ref{320e}) and Corollary  \ref{21c} $(i)$  it follows that
\begin{theorem}\label{41t} Let $\pi$ and $\pi'$ be two different
non-increasing tree degree sequences with $\pi\triangleleft \pi'$.
 Let $T$ and $T'$ be an extremal tree of  $\Gamma(\pi)$ and $\Gamma(\pi')$, respectively.   If  $f(x,y)$ is a good  escalating function, then $M_{f}(T)< M_{f}(T').$
\end{theorem}

 \begin{theorem}\label{42t} Let $\pi$ and $\pi'$ be two different
non-increasing unicyclic  degree sequences with $\pi\triangleleft \pi'$.
 Let $U$ and $U'$ be an extremal unicyclic graph of  $\Gamma(\pi)$ and $\Gamma(\pi')$, respectively.   If  $f(x,y)$ is a good  escalating function, then $M_{f}(U)< M_{f}(U').$
\end{theorem}
  {\bf Proof.}   Suppose that  $\pi=(d_{1},d_{2},...,d_{n})$ and $\pi'=(d'_{1},d'_{2},...,d'_{n})$. Since   $\pi\lhd \pi'$, by Lemma \ref{26l} we may suppose that $\pi$ and $\pi'$ differ
only in two positions, where the differences are  $1$. Thus,  we may assume that
$d_{i}=d '_{i}$ for $i\neq p,q$, and $d_{p}+1=d'_{p}$, $d_{q}-1=d'_{q}$.  \par

 If $d_{n}=2$, then $\pi=(2,2,...,2)$ and $\pi'=(3,2,...,2,1)$. By Corollary \ref{21c} $(ii)$, it suffices to show that  $M_{f}(C_{n})< M_{f}(U_{M}(\pi')).$  Since $U_{M}(\pi')$ contains at least one  pendant  vertex, $n\geq 4$. We may suppose that $n\geq 5$, as the case of $n=4$ can be proved similarly. Let $g(x,y)=\frac{\partial f(x,y)}{\partial x}$. In this case, since $\frac{\partial^{2} f(x,y)}{\partial x^{2}}\geq 0$ and $g(x,y)>0$,  by Corollary \ref{21c} $(ii)$ we have   \begin{align*}  M_{f}(U_{M}(\pi'))-M_{f}(C_{n})=&3f\left(3,2\right)+f(2,1)-4f(2,2)\\>&f\left(3,2\right)+f(2,1)-2f(2,2)\\=&f\left(3,2\right)-f(2,2)-(f(2,2)-f(1,2))
 \\=&g(2+a,2)-g(1+b,2)\geq 0,\end{align*}
where $0\leq a\leq 1$ and $0\leq b\leq 1$.

 Otherwise,   $d_{n}=1$. In this case,  it suffices to show that $M_{f}(U_{M}(\pi))< M_{f}(U_{M}(\pi'))$ by    Corollary \ref{21c} $(ii)$. If $2\leq q\leq 3$, then $d_{q}\geq 3$ since $U_{M}(\pi')\in \Gamma(\pi')$. If $q\geq 4$, then $d_{q}\geq 2$.  In both cases,    $U_{M}(\pi)$ contains a surprising vertex. Now, the result follows from (\ref{320e}).   \qed

\begin{theorem}\label{43t} Let $\pi$ and $\pi'$ be two different
non-increasing bicyclic  degree sequences with $\pi\triangleleft \pi'$.
 Let $B$ and $B'$ be an extremal bicyclic  graph of  $\Gamma(\pi)$ and $\Gamma(\pi')$, respectively.   If  $f(x,y)$ is a good  escalating function, then $M_{f}(B)< M_{f}(B').$
\end{theorem}
  {\bf Proof.}   Suppose that  $\pi=(d_{1},d_{2},...,d_{n})$ and $\pi'=(d'_{1},d'_{2},...,d'_{n})$. Since   $\pi\lhd \pi'$, by Lemma \ref{26l} we may suppose that $\pi$ and $\pi'$ differ
only in two positions, where the differences are  $1$. Thus,  we may assume that
$d_{i}=d '_{i}$ for $i\neq p,q$, and $d_{p}+1=d'_{p}$, $d_{q}-1=d'_{q}$. It suffices to show that $M_{f}(B_{M}(\pi))< M_{f}(B_{M}(\pi'))$ by  Theorem \ref{22t}. \par

We firstly consider the case of  $d_{n}\geq 2$. In this case,  either $\pi=(4,2^{(n-1)})$ or $\pi=(3^{(2)},2^{(n-2)})$. If $\pi=(4,2^{(n-1)})$, then either $\pi'=(5,2^{(n-2)},1)$ or $\pi'=(4,3,2^{(n-3)},1)$. If $\pi=(3^{(2)},2^{(n-2)})$, then $\pi'=(4,2^{(n-1)})$ or $\pi'=(3^{(3)},2^{(n-4)},1)$ or $\pi'=(4,3,2^{(n-3)},1)$. Whenever which case happens, according to Theorem  \ref{22t}, $B_{M}(\pi)$ contains a surprising vertex, and hence the result follows from (\ref{320e}). We secondly consider the  case of  $d_{n}=1$.  \par
{\bf Case  1.}   $d_{n}=1$ and  $d_{2}=2$. \par

In this case, $\pi=(k+4,2^{(n-k-1)},1^{(k)})$, where $k\geq 1$. Then, either $\pi'=(k+5,2^{(n-k-2)},1^{(k+1)})$ or $\pi'=(k+4,3,2^{(n-k-3)},1^{(k+1)})$. According to Theorem \ref{22t},  $\mathcal {R}(B_{M}(\pi))\cong B_2$. \par

 If $\pi'=(k+4,3,2^{(n-k-3)},1^{(k+1)})$, since $v_{5}\in N_{B_{M}(\pi)}(v_{4})\setminus N_{B_{M}(\pi)}[v_{2}]$ and  $v_{5}\not\in V(P_{v_{2}v_{4}})$, it is easy to see that  $v_{5}$ is a surprising vertex of  $B_{M}(\pi)$.

If $\pi'=(k+5,2^{(n-k-2)},1^{(k+1)})$, then $p=1$,  $q\geq 6$ and $d_{q}=2$. By the structure of $B_{M}(\pi)$, there exists a vertex   $v_{k}$ $(k>q)$ such that $v_{k}\in N_{B_{M}(\pi)}(v_{q})\setminus N_{B_{M}(\pi)}[v_{p}]$ and  $v_{k}\not\in V(P_{v_{p}v_{q}})$. In this case, it is easy to see that  $v_{k}$ is a surprising vertex of $B_{M}(\pi)$.

In both cases, the result follows from (\ref{320e}).
 \par
{\bf Case  2.} $d_{n} = 1$ and  $d_{1}\geq  d_{2}\geq 3$.\par

By Theorem \ref{22t}, $\mathcal {R}(B_{M}(\pi))\cong B_1$. If $ q\geq5$, then $d_{q}=d'_{q}+1\geq 2$. If $3\leq q\leq 4$, since $B_{M}(\pi')$ is a bicyclic graph, $d_{q}=d'_{q}+1\geq 3$. In both cases,  by the structure of $B_{M}(\pi)$, there exists vertex $v_{k}$ $(k>q)$ such that $v_{k}\in N_{B_{M}(\pi)}(v_{q})\setminus  N_{B_{M}(\pi)}[v_{p}]$ and $v_{k}\not\in V(P_{v_{p}v_{q}})$. In this case, it is easy to see that  $v_{k}$ is a surprising vertex of  $B_{M}(\pi)$. Now, the result follows from (\ref{320e}). Thus, we may suppose that $q=2$, and hence $p=1$ in the sequel.
\par
If $q=2$ and $d_{2}\geq 4$, since $\mathcal {R}(B_{M}(\pi))\cong B_1$,  there exists a vertex $v_{k}\in N_{B_{M}(\pi)}(v_{2})\setminus N_{B_{M}(\pi)}[v_{1}]$ and $v_{k}\not\in V(P_{v_{1}v_{2}})$,  which implies that $v_{k}$ is a surprising vertex of $B_{M}(\pi)$. Now, the result follows from (\ref{320e}). \par
If $q=2$ and $d_{2}=3$, then $d'_{2}=2$. Thus,  $\pi=(k+3,3,2^{(n-k-2)},1^{(k)})$ and $ \pi'=(k+4,2^{(n-k-1)},1^{(k)})$, where $k\geq 1$. By Theorem \ref{22t}, $\mathcal {R}(B_{M}(\pi))\cong B_1$ and $\mathcal {R}(B_{M}(\pi'))\cong B_2$.

Then,  $B_{M}(\pi)$ is obtained from $B_1$ by attaching  $k$ paths, say $P_{l_1}$, $P_{l_2}$, $...,$ $P_{l_k}$,  of almost equal lengths to one vertex of degree three  of $B_1$, where $l_1\geq l_2\geq \cdots\geq l_k\geq 1$ and $l_1\geq 2$.

 If $l_1\geq 3$, since $g(x,y)=\frac{\partial f(x,y)}{\partial x}>0$ and
 $\frac{\partial^{2} f(x,y)}{\partial x^{2}}\geq 0$, we have  $M_f(B_{M}(\pi'))-M_f(B_{M}(\pi))=(k+4)f(k+4,2)+f(2,2)-(k+2)f(k+3,2)-f(k+3,3)-2f(3,2)=(k+2)g(k+3+a,2)-2f(3,2)+2f(k+4,2)+f(2,2)-f(k+3,3)
=(k+2)g(k+3+a,2)-g(2+b,2)+2f(k+4,2)-f(k+3,3)-f(3,2)>2f(k+4,2)-f(k+3,3)-f(3,2)\geq 0,$ where $0\leq a\leq 1$ and $0\leq b\leq 1$.

Otherwise, $l_1=2.$ We may suppose that  $(l_1,l_2,..., l_k)=(2^{(s)},1^{(k-s)})$, where $1\leq s\leq k$. In this case,
$M_f(B_{M}(\pi'))-M_f(B_{M}(\pi))=(s+3)f(k+4,2)+2f(2,2)+(k+1-s)f(k+4,1)-(s+2)f(k+3,2)-f(k+3,3)-2f(3,2)-(k-s)f(k+3,1)-f(2,1)
>(s+2)g(k+3+a,2)-2g(2+b,2)+f(k+4,2)+f(k+4,1)-f(k+3,3)-f(2,1)
>f(k+4,2)+f(k+4,1)-f(k+3,3)-f(2,1)\geq f(k+3,3)+f(3,1)-f(k+3,3)-f(2,1)>0,$ where $0\leq a\leq 1$ and $0\leq b\leq 1$.
\qed

\par\bigskip
Let $F_{n}(k)$ (respectively, $F'_{n}(k)$, $F''_{n}(k)$) be the tree (respectively,
unicyclic graph, bicyclic graph)  on $n$ vertices obtained  by attaching $k$ paths of almost
equal lengths to one isolated vertex (respectively, one vertex of  $C_{3}$, the vertex of degree four of $B_2$).  From Theorems \ref{21t}--\ref{22t} and  Theorems \ref{41t}--\ref{43t}, we can deduce the following extremal results easily:
 \begin{corollary} \label{31c} If $k \geq 1$ and $f(x,y)$ is   good  escalating, then
    $F''_{n}(k)$ $($respectively, $F'_{n}(k)$, $F_{n}(k)$$)$  is an extremal  bicyclic graph $($respectively, unicyclic graph, tree$)$  among all bicyclic graphs  $($respectively, unicyclic graphs, trees$)$ with  $n$ vertices and  $k$ pendant vertices when $n\geq k+5$  $($respectively, $n\geq k+3$, $n\geq k+1$$)$.
\end{corollary}
{\bf Proof.} Here, we only prove the case of bicyclic graphs, as the cases of unicyclic graphs and trees  can be proved similarly.
 Let $\pi'=(k+4,2^{(n-k-1)},1^{(k)})$. Suppose that $G$ is a bicyclic graph with $\pi=(d_{1},d_{2},...,d_{n})$ as its degree sequence, where $d_{n}=d_{n-1}=\cdots=d_{n-k+1}=1$ and $d_{n-k}\geq 2$. If $\pi\neq \pi'$, then
$\pi\lhd \pi'$. Now, the result  follows  from  Theorem  \ref{22t} and Theorem \ref{43t}. \qed
\par \bigskip

Let $S^{(0)}_n$ be the star with $n$ vertices, let  $S^{(1)}_n$   be the unicyclic graph  obtained from the star  $S^{(0)}_n$   by adding one edge  between two pendant vertices of $S^{(0)}_n$, and let $S^{(2)}_n$   be the bicyclic graph  obtained from  $B_1$   by attaching $n-4$ pendant vertices to one vertex of degree three of   $B_1$.
 \begin{corollary} \label{32c} If $n\geq 4$ and $f(x,y)$ is   good  escalating, then   $S^{(c)}_n$  is the unique  extremal  $c$-cyclic graph  among all $c$-cyclic graphs   with  $n$ vertices for $0\leq c\leq 2$.
\end{corollary}
{\bf Proof.} Here, we only prove the case of bicyclic graphs, as the cases of unicyclic graphs and trees  can be proved similarly.  Let   $\pi'=(n-1,3,2^{(2)},1^{(n-4)})$. Then, $\pi'$   uniquely maximizes those degree sequences of bicyclic graphs with $n$ vertices in the relation $\lhd$. Furthermore,    $S^{(2)}_{n}$ is the unique bicyclic graphs of  $\Gamma(\pi')$. Now, the result  follows  from Theorem \ref{43t}.   \qed
\section{Some applications}
  As an extension of $Z_2(G)$, we define $Z_\alpha(G)$ as follows:
\begin{align}\label{41e}Z_{\alpha}(G)=\sum_{uv\in E(G)}\left(d(u)+d(v)-2\right)^{\alpha}.\end{align}

In  \cite{Ba1}, Zhang et al.  has shown that $\chi_{\alpha}(G)$ is    escalating    for   $\alpha\geq 1$ and  $\chi_{\alpha}(G)$ is  de-escalating for $0<\alpha<1$, and  $W_{\alpha}(G)$ is    escalating    for   $\alpha>0$. In the following, we will  extend  these results.

 \begin{theorem}\label{31t} Let $G$ be a connected graph with at least three vertices. Then,    $(i)$ $Z_{\alpha}(G)$ and $\chi_{\alpha}(G)$ are   escalating  for \textcolor{red}{$\alpha>1$} and $\alpha<0$. Furthermore,  $Z_{\alpha}(G)$ and $\chi_{\alpha}(G)$ are  good   escalating    for \textcolor{red}{$\alpha>1$}. $(ii)$ $Z_{\alpha}(G)$ and $\chi_{\alpha}(G)$ are   de-escalating for $0<\alpha<1$.
\end{theorem}
{\bf Proof.} Here, we only   prove  $(i)$ and $Z_{\alpha}(G)$, as $(ii)$ and $\chi_{\alpha}(G)$ can be proved similarly. By (\ref{41e}), we may define   $f(x, y)=(x+y-2)^{\alpha}$ with  $\min\{x,y\}\geq 1$ and $\max\{x,y\}\geq 2$. Suppose that  $x_1 \geq y_1$ and $x_2 \geq y_2$. To show that $f(x, y)=(x+y-2)^{\alpha}$ is escalating, by (\ref{21e}) it suffices to show that $$\int_{y_1+x_2-2}^{x_1+x_2-2} \alpha t^{\alpha-1} dt\geq \int_{y_1+y_2-2}^{x_1+y_2-2} \alpha t^{\alpha-1} dt,$$ which is equivalent to \begin{align}\label{31e}\int_{y_1}^{x_1} \alpha (t+x_2-2)^{\alpha-1} dt\geq \int_{y_1}^{x_1} \alpha (t+y_2-2)^{\alpha-1} dt.\end{align} When \textcolor{red}{$\alpha >1$} or $\alpha<0$, since  $t+x_2-2\geq t+y_2-2\geq y_1+y_2-2\geq 1$, we can conclude that inequality (\ref{31e}) holds.
Thus, $f(x,y)$ is escalating.  \par When \textcolor{red}{$\alpha >1$}, it is easy to see that  $f(x_1+1,x_2)+f(x_1+1,y_1-1)-f(x_2,y_1)-f(x_1,y_1)=(x_1+x_2-1)^{\alpha}-(y_1+x_2-2)^{\alpha}\geq0$   holds for any $x_1\geq y_1$  and $x_2\geq 1$ (since $\max\{x_1,y_1\}\geq 2$, $\max\{x_2,y_1\}\geq 2$). Furthermore, since  $x+y\geq 3$,    $ g(x,y)=\frac{\partial f(x,y)}{\partial x}=\alpha(x+y-2)^{\alpha-1}>0$ and $\frac{\partial^{2} f(x,y)}{\partial x^{2}}=\alpha(\alpha-1)(x+y-2)^{\alpha-2}\geq 0$.   Thus,  $Z_{\alpha}(G)$ is a  good  escalating function for \textcolor{red}{$\alpha>1$}. \qed

\begin{theorem}\label{33t} Let $G$ be a connected graph with at least three vertices. Then, $W_{\alpha}(G)$ is    escalating    for any real number $\alpha\neq 0$. Furthermore, $W_{\alpha}(G)$ is  good  escalating for $\alpha=1$.
\end{theorem}
{\bf Proof.} Let $f(x,y)=(xy)^{\alpha}$ and $g(x,y)=\frac{\partial f(x,y)}{\partial x}$, where $\alpha\neq 0$, $\min\{x,y\}\geq 1$ and $\max\{x,y\}\geq 2$. For any real number $\alpha$,  it is easy to see that $$f(x_1,x_2)+f(y_1,y_2)-f(x_1,y_2)-f(x_2,y_1)=\left(x^{\alpha}_1-y^{\alpha}_1\right)\left(x^{\alpha}_2-y^{\alpha}_2\right)\geq 0$$
holds for any $x_1 \geq y_1 \geq 1$ and $x_2 \geq y_2\geq 1$, and hence $W_{\alpha}(G)$ is    escalating. \par
It is easy to check that $W_{\alpha}(G)$ is  good  escalating for $\alpha=1$.  \qed

\begin{remark}\label{41c} {\em Note that  $W_{1}(G)=M_2(G)$. By Theorems \ref{41t}--\ref{43t} and Theorem \ref{33t}, we can easily   deduce  the main results of \cite{Mh1,Mh2,Zhang2}, that is  ``If $\pi\lhd \pi'$, then    $M_{2}(G')> M_{2}(G)$  holds for any two extremal graphs $G\in\Gamma(\pi)$ and $G'\in\Gamma(\pi')$ when $G$ is a $c$-cyclic graph for $0\leq c\leq 2$."}
\end{remark}

\begin{remark}\label{42c} {\em   By Theorems \ref{21c} and \ref{33t}, we can conclude that an  extremal graph $G$ of the Randi\'c index  or second Zagreb index  among these connected graphs with fixed degree sequence   is a $BFS$-tree \cite{Wang2,Mh2} (respectively, $U_M(\pi)$, $B_M(\pi)$) when $G$ is   a tree (respectively, a unicyclic graph, a bicyclic graph  with $d_n=1$ \cite{Mh1,Zhang2}).}
\end{remark}

\begin{remark}\label{43c} {\em   By Theorems \ref{21c} and \ref{22t}, we can conclude that an  extremal graph $G$ of sum-connectivity index, the third Zagreb  index, reformulated Zagreb index or harmonic index  among these connected graphs with fixed degree sequence   is a $BFS$-tree (respectively, $U_M(\pi)$, $B_M(\pi)$) when $G$ is   a tree (respectively, a unicyclic graph, a bicyclic graph  with $d_n=1$). Furthermore, Theorem \ref{31t} implies that Theorems \ref{41t}--\ref{43t} are also  suitable for the  third Zagreb  index or reformulated Zagreb index.}
\end{remark}

\begin{remark}\label{43c} {\em Zhang et al. \cite{Ba1} has shown that  $f(x,y)=\sqrt{\frac{x+y-2}{xy}}$ is de-escalating, and hence the results of Theorems  \ref{21c}--\ref{22t} are also  suitable for the Atom-Bond connectivity   index   of $G$ \cite{Zhou3,Lin1}. }
\end{remark}
\begin{remark}\label{44c} {\em By Theorem \ref{31t}, $Z_2(G)$ is a good escalating function. Thus, Corollary \ref{32c} implies that   $S^{(c)}_n$  is the unique  extremal  $c$-cyclic graph with maximum reformulated Zagreb indices  among all $c$-cyclic graphs   with  $n\geq 4$ vertices for $0\leq c\leq 2$. These results had been proved in \cite{I1,J1}. }
\end{remark}

\noindent
 {\it Acknowledgement\/}. The authors would like  to thank  the referees
for their valuable comments which lead to an improvement of the
original manuscript.

\end{document}